\documentclass[11pt]{amsart}

\usepackage{graphicx}
\numberwithin{equation}{section}
\numberwithin{figure}{section}
\theoremstyle{definition}
\newtheorem{remark}{Remark}
\numberwithin{remark}{section}

\newcommand{\lsp}{\vspace{3mm}}

\begin{document}

\begin{center}
\textbf{\large A fast direct solver for network matrices}

\lsp

P.G. Martinsson, Dept.~of Applied Mathematics, University of Colorado at Boulder

\lsp

\begin{minipage}{125mm}
\textbf{Abstract:} A fast direct inversion scheme
for the large sparse systems of linear equations resulting from the discretization
of elliptic partial differential equations in two dimensions is given. The scheme is
described for the particular case of a discretization on a uniform square
grid, but can be generalized to more general geometries. For a grid
containing $N$ points, the scheme requires $O(N \log^{2}N)$ arithmetic operations
and $O(N \log N)$ storage to compute an approximate inverse. If only a single
solve is required, then the scheme requires only $O(\sqrt{N}\,\log N)$ storage;
the same storage is sufficient for computing the Dirichlet-to-Neumann operator
as well as other boundary-to-boundary operators.
The scheme is illustrated with several numerical examples. For instance,
a matrix of size $10^6 \times 10^6$ is inverted to seven digits accuracy
in four minutes on a 2.8GHz P4 desktop PC.

\end{minipage}

\end{center}

\section{Introduction}

This note describes a scheme for rapidly solving the systems of linear
equations arising from the finite element or finite difference discretization
of elliptic partial differential equations in two dimensions. Unlike most
existing fast schemes
which rely on iterative solvers (GMRES / conjugate gradient / \dots), the
method described here directly inverts the matrix of the linear system. This
obviates the need for customized pre-conditioners, and leads to dramatic
speed-ups in environments in which the same equation needs to be solved for
a sequence of different right-hand sides.

The scheme is described for the case of equations defined on a uniform square
grid. Slight modifications would enable the handling of uniform grids on fairly
regular two-dimensional domains, but further work would be required to construct
methods for non-uniform grids or complicated geometries.

For a system matrix of size $N\times N$ (corresponding to a
$\sqrt{N} \times \sqrt{N}$ grid), the scheme requires
$O(N\log^{2}N)$ arithmetic operations, and $O(N\log N)$ storage to
compute the full inverse. For a single solve, only $O(\sqrt{N}\,\log
N)$ storage is required. Moreover, for problems loaded on the
boundary only, any solves beyond the first require only
$O(\sqrt{N}\,\log N)$ arithmetic operations (provided that only the
solution on the boundary is sought). Numerical experiments indicate
that the constants in these asymptotic estimates are quite moderate.
For instance, to directly solve a system involving a $1\,000\,000
\times 1\,000\,000$ matrix to seven digits of accuracy takes about
four minutes on a 2.8GHz desktop PC with 512Mb of memory. Additional
solves beyond the first can be performed in 0.03 seconds (provided
only boundary data is involved).

The scheme is inspired by earlier work by Hackbusch and co-workers
\cite{hackbusch,H2_matrix,hackbusch2003,bebendorf2003} and
work by Gu and Chandrasekaran \cite{chan_gu_2001,chan_gu_2003,gu_SSS,HSS_sparse}.
These authors have published methods that
rapidly perform algebraic operations (matrix-vector multiplies, matrix-matrix multiplies,
matrix inversion, \textit{etc}) on matrices whose off-diagonal blocks
have low rank. Our scheme is similar, but uses simpler data structures.

The fast inversion scheme retains its $O(N\log^{2}N)$ computational complexity for
a wide range of network matrices. It does not rely on the fact that the matrix is
associated with a PDE; rather, it works for any matrix whose inverse has off-diagonal
blocks of low rank (what Hackbusch calls an $\mathcal{H}$-matrix, and Gu and Chandrasekaran
calls an ``HSS-matrix'').

The paper is structured as follows: Section \ref{sec:prel} describes
some known algorithms for performing fast operations on matrices with
off-diagonal blocks of low rank. Section \ref{sec:babuska} describes
a very simple $O(N^{2})$ inversion scheme. Section \ref{sec:fast}
describes how the $O(N^{2})$ scheme of Section \ref{sec:babuska} can be
accelerated to $O(N\log^{2}N)$ using the methods of Section \ref{sec:prel}.
Section \ref{sec:numerics} give the results of numerical experiments.

\section{Preliminaries}
\label{sec:prel}

In this section, we summarily describe a class of matrices for which
there are known algorithms for rapidly evaluating the result of
algebraic operations such as matrix-vector products, matrix-matrix
products, and matrix inversions. The basic concept is that the
off-diagonal blocks of the matrix can, to within some preset
accuracy $\varepsilon$, be approximated by low-rank matrices.
Different authors have given different definitions, and provided
different algorithms; we mention that Hackbusch refers to such
matrices as $\mathcal{H}$-matrices \cite{hackbusch}, or
$\mathcal{H}^{2}$-matrices \cite{H2_matrix}, while Gu and Chandrasekaran refers to
them as ``Hierarchically Semi-Separable'' (HSS) matrices \cite{HSS_sparse} or
``Sequentially Semi-Separable'' (SSS) matrices \cite{gu_SSS}.

\begin{figure}
\includegraphics[height=80mm]{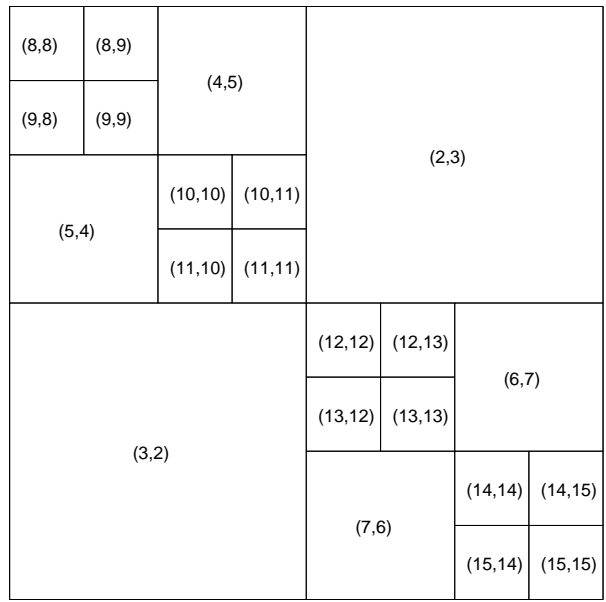}
\caption{Tessellation of an HSS matrix.}
\label{fig:tessellation}
\end{figure}

In this paper we will use the term HSS matrix to refer to an
$n\times n$ matrix $B$ that can be tessellated in the pattern shown
in Figure \ref{fig:tessellation} in such a way that the rank of every block in the
tessellation is bounded by some fixed small number $p$. Such a
matrix can clearly be stored in $O(p\,n\log n)$ storage, and can be
applied to a vector in $O(p\,n \log n)$ arithmetic operations.
Moreover, via a trivial recursion (described in Appendix
\ref{sec:invert}), it is possible to invert such a matrix in
$O(p\,n\log^{2}n)$ operations.

\begin{remark}
We will in this paper not distinguish between a matrix that is
exactly an HSS-matrix, and a matrix that can to high precision be
approximated by an HSS-matrix.
\end{remark}

\begin{remark}
\label{remark:SSS}
We use the term SSS-matrix to refer to a matrix that is an
HSS-matrix but has the additional property the bases for the row and
column spaces of the off-diagonal blocks can be expressed
hierarchically. As an example, a basis for the column space of the
block labeled $(4,5)$ is constructed from the bases for the column
spaces of the blocks $(8,9)$ and $(9,8)$. An SSS-matrix requires
only $O(N)$ storage, and can be applied to a vector in $O(N)$
operations. Matrix-inversion schemes for SSS-matrices are a big more
complicated than HSS inversion schemes, but in some environments,
$O(N)$ inversion schemes exist, see \textit{e.g.}~\cite{m2003}.
\end{remark}

\begin{remark}
In Hackbusch's terminology, what we call an HSS-matrix roughly
corresponds to an $\mathcal{H}$-matrix, and what we call an
SSS-matrix, roughly corresponds to an $\mathcal{H}^{2}$-matrix. The
tessellations used by Hackbusch are slightly different, though.
\end{remark}

\section{An exact $O(N^{2})$ inversion scheme}
\label{sec:babuska}

\begin{figure}
\includegraphics[height=100mm]{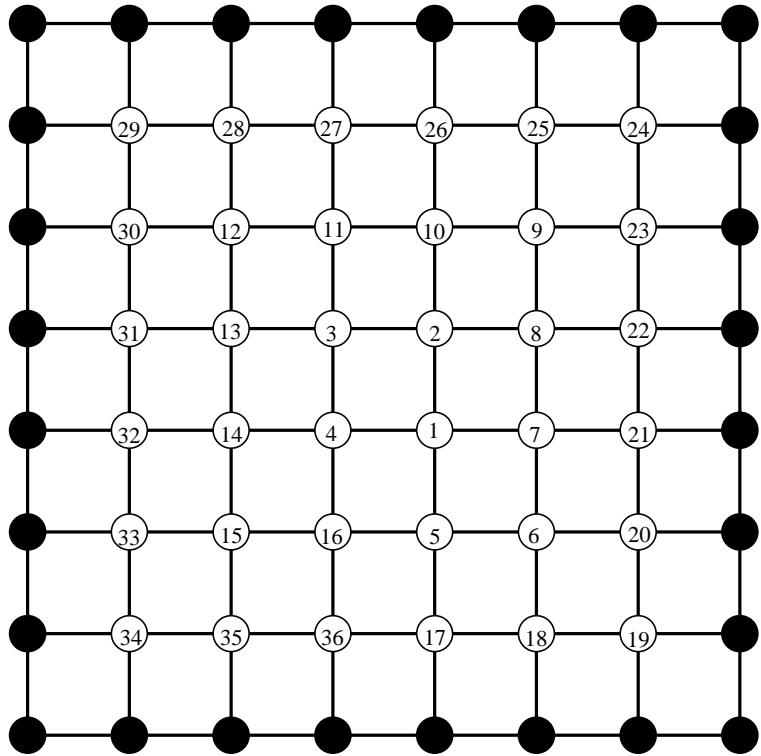}
\caption{The computational grid for $n = 6$.}
\label{fig:grid}
\end{figure}

Letting $n$ be a positive even integer, we consider a static conduction
problem on a square $(n+2)\times (n+2)$ grid such as the one illustrated
in Figure \ref{fig:grid}. Each interior node is connected to its four nearest neighbors
with bars with positive conductivity. We prescribe the temperatures on the
boundary nodes (marked with filled circles in Figure \ref{fig:grid}), and seek the equilibrium
temperatures of the $N$ internal nodes, where $N = n^{2}$. This problem can
be formulated as a linear system
\begin{equation}
\label{eq:Ax=b}
A\,x = b,
\end{equation}
where $A$ is an $N\times N$ sparse matrix, $x$ is an $N\times 1$ vector of
unknown temperatures, and $b$ is an $N\times 1$ vector derived from the
pre-scribed boundary temperatures. When all bars have unit conductivity,
the matrix $A$ in (\ref{eq:Ax=b}) is the well-known five-point stencil
with coefficients $4$/$-1$/$-1$/$-1$/$-1$.

In this section we describe a method for directly solving the linear
system (\ref{eq:Ax=b}) that relies on the sparsity structure of
the matrix only. In the absence of rounding
errors, it would be exact. When the matrix $A$ is of size $N\times N$, the
scheme of this section requires $O(N^{2})$ floating point operations and
$O(N)$ memory. This makes the scheme significantly slower
than well-known $O(N^{3/2})$ schemes such as nested dissection.
(We mention that $O(N^{3/2})$ is optimal in
this environment.) The only purpose of the scheme presented
in this section is that it can straight-forwardly be accelerated to
an $O(N)$ or $O(N \log^{\kappa} N)$ scheme, as shown in Section \ref{sec:fast}.

\begin{figure}
\includegraphics[height=100mm]{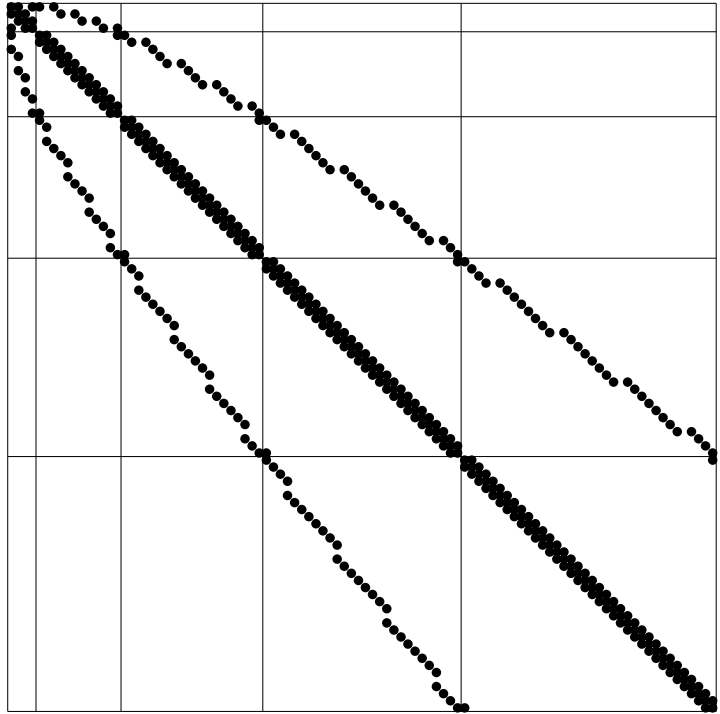}
\caption{The sparsity pattern of $A$ in (\ref{eq:Ax=b}) for $N = 100$.}
\label{fig:sparsity}
\end{figure}

Ordering the $N$ points in the grid in the spiral pattern shown in Figure \ref{fig:grid},
the matrix $A$ in equation (\ref{eq:Ax=b}) has the sparsity pattern shown in
Figure \ref{fig:sparsity}. We next partition the grid into $n$ concentric squares, and collect
the nodes into index sets $J_{1},\, J_{2},\,\dots\,,J_{n}$ accordingly. In other words,
\begin{equation*}
\begin{array}{rcl}
J_{1} &=& \{1,\,2,\,3,\,4\},\\
J_{2} &=& \{5,\,6,\,\dots\,,16\},\\
&\vdots & \\
J_{n} &=& \{(2n-1)^{2} + 1,\,(2n-1)^{2} + 2,\,\dots\,,(2n)^{2}\}.
\end{array}
\end{equation*}
For $\kappa,\lambda \in \{1,\,2,\,\dots\,n\}$, we let $A_{\kappa\lambda}$ denote the submatrix of
$A$ formed by the intersection of the $J_{\kappa}$ rows with the $J_{\lambda}$ columns.
Then the linear system (\ref{eq:Ax=b}) takes on the block-tridiagonal form (also shown in Figure \ref{fig:sparsity})
\begin{equation}
\label{eq:Ax=b_blocked}
\left[\begin{array}{cccccc}
A_{11} & A_{12} & 0      & 0      & \cdots  & 0\\
A_{21} & A_{22} & A_{23} & 0      & \cdots  & 0\\
0      & A_{32} & A_{33} & A_{34} & \cdots  & 0\\
0      & 0      & A_{43} & A_{44} & \cdots  & 0\\
\vdots & \vdots & \vdots & \vdots &         & \vdots\\
0      & 0      & 0      & 0      & \cdots  & A_{nn}
\end{array}\right]\,
\left[\begin{array}{c}
x_{1} \\ x_{2} \\ x_{3} \\ x_{4} \\ \vdots \\ x_{n}
\end{array}\right] =
\left[\begin{array}{c}
b_{1} \\ b_{2} \\ b_{3} \\ b_{4} \\ \vdots \\ b_{n}
\end{array}\right],
\end{equation}
where $x$ and $b$ have been split accordingly.

The equation (\ref{eq:Ax=b_blocked}) can now easily be solved by
eliminating the variables one by one. Using the first equation to eliminate
$x_{1}$ from the second one, we obtain the following system of
equations for the variables $x_{2},\dots,x_{n}$:
\begin{equation}
\left[\begin{array}{cccccc}
\tilde{A}_{22} & A_{23} & 0      & \cdots  & 0\\
A_{32}         & A_{33} & A_{34} & \cdots  & 0\\
0              & A_{43} & A_{44} & \cdots  & 0\\
\vdots         & \vdots & \vdots &         & \vdots\\
0              & 0      & 0      & \cdots  & A_{nn}
\end{array}\right]\,
\left[\begin{array}{c}
x_{2} \\ x_{3} \\ x_{4} \\ \vdots \\ x_{n}
\end{array}\right] =
\left[\begin{array}{c}
\tilde{b}_{2} \\ b_{3} \\ b_{4} \\ \vdots \\ b_{n}
\end{array}\right],
\end{equation}
where $\tilde{A}_{22} = A_{22} - A_{21}A_{11}^{-1}A_{12}$ and
$\tilde{b}_{2} = b_{2} - A_{21}A_{11}^{-1}b_{1}$. The process
used to eliminate $x_{1}$ can easily be continued to eliminate
the first $n-1$ blocks. This leaves us with the
$(8n-4)\times(8n-4)$ system
$$
\tilde{A}_{nn}\,x_{n} = \tilde{b}_{n},
$$
which we solve directly to obtain $x_{n}$. Once $x_{n}$ is known, we compute
$x_{n-1}$ by solving the system
$$
\tilde{A}_{n-1,n-1}\,x_{n-1} = \tilde{b}_{n-1} - A_{n-1,n}\,x_{n}.
$$
The remaining $x_{j}$'s are computed analogously. The entire process
is summed up in Algorithm I.

\begin{figure}
\begin{center}
\fbox{
\begin{minipage}{120mm}
\begin{center}
\underline{\textit{Algorithm 1:}}
\end{center}

\lsp

\begin{tabbing}
\hspace{10mm} \= \hspace{3mm} \= \kill
(1) \> $\tilde{A}_{11} = A_{11}$ and $\tilde{b}_{1} = b_{1}$.\\
\\
(2) \> \textbf{for} $\kappa = 2:n$\\
(3) \>\> $\tilde{A}_{\kappa\kappa} = A_{\kappa\kappa} -
A_{\kappa,\kappa-1}\,\tilde{A}_{\kappa-1,\kappa-1}^{-1}\,A_{\kappa-1,\kappa}$\\
(4)\>\>$\tilde{b}_{\kappa} = b_{\kappa} -
A_{\kappa,\kappa-1}\,\tilde{A}_{\kappa-1,\kappa-1}^{-1}\,\tilde{b}_{\kappa-1}$\\
(5)\>\textbf{end}\\
\\
(6)\>$x_{n} = \tilde{A}_{nn}^{-1}\,\tilde{b}_{n}$\\
\\
(7)\>\textbf{for} $\kappa = (n-1):(-1):1$\\
(8)\>\>$x_{\kappa} = \tilde{A}_{\kappa\kappa}^{-1}\,\bigl(\tilde{b}_{\kappa} - A_{\kappa,\kappa+1}\,x_{\kappa+1}\bigr)$\\
(9)\>\textbf{end}
\end{tabbing}
\end{minipage}}
\end{center}
\caption{This algorithm directly solves the tridiagonal system
of equations (\ref{eq:Ax=b_blocked}).}
\label{fig:algorithm1}
\end{figure}

We note that while all matrices $A_{\kappa\lambda}$ are sparse,
the matrices $\tilde{A}_{\kappa\kappa}$ are dense. This means that the
cost of inverting $\tilde{A}_{\kappa\kappa}$ in each step of Algorithm I
is $O(\kappa^{3})$. (Note that the remaining matrix-matrix operations
involve matrices that are diagonal or tri-diagonal and have negligible
costs in comparison to the matrix inversion.) The total cost $T_{\rm total}^{\rm dense}$
is therefore
$$
T_{\rm total}^{\rm dense} \sim \sum_{\kappa=1}^{n}\kappa^{3} \sim n^{4} \sim N^{2}.
$$

\begin{remark}
\label{remark:cheap}
We note that the matrices $\tilde{A}_{\kappa\kappa}$ are not merely artificial
objects of the numerical algorithm. In fact, $\tilde{A}_{nn}^{-1}$ is the matrix
that maps the load on the boundary to the potential on the boundary. In many
environments, the interior of the grid is unloaded (\textit{i.e.}~$b_{\kappa} = 0$
when $\kappa = 1,\,2,\dots,n-1$). In such cases, the operator $\tilde{A}_{nn}^{-1}$
is the solution operator of the original problem. Since none of the intermediate
matrices $\tilde{A}_{\kappa\kappa}$ are required in this environment, the
algorithm only requires $O(N)$ memory.
\end{remark}

\section{A fast direct solver}
\label{sec:fast}

The computational cost of Algorithm I consists almost entirely of
the inversion of the dense matrices $\tilde{A}_{\kappa\kappa}$.
This cost can be greatly reduced whenever $\tilde{A}_{\kappa\kappa}$
is a ``compressible'' matrix in the sense described in Section
\ref{sec:prel}. It is not the purpose of this paper to classify
exactly when this happens; we will simply note that $A$ is
typically compressible in this sense when it results from the
discretization of an elliptic PDE, while it is typically not compressible
when it results from the discretization of a highly oscillatory wave equation.
We also note that the condition that $A$
result from the discretization of an elliptic PDE is very far from necessary,
as the numerical examples in Section \ref{sec:numerics} demonstrate.

In cases where each $\tilde{A}_{\kappa\kappa}$ is what we in Section \ref{sec:prel}
called in ``HSS''-matrix,
very simple inversion schemes (see Appendix \ref{sec:invert}) are available
for computing an approximation to $\tilde{A}_{\kappa\kappa}^{-1}$ in $O(\kappa\,\log^{2}\kappa)$
arithmetic operations. The total computational cost of Algorithm I is then
$$
T_{\rm total}^{\rm HSS} \sim \sum_{\kappa=1}^{n} \kappa\,(\log\kappa)^{2} \sim
n^{2}\,(\log n)^{2} \sim N\,(\log N)^{2}.
$$
This scheme requires $O(N\,\log N)$ memory to store the entire approximation to $A^{-1}$.
Such a scheme has been implemented and the computational results are given in Section \ref{sec:numerics}.

We note that the scheme is particularly memory efficient in environments where
the domain is loaded on the boundary only, and only the solution at the boundary
is sought, \textit{cf.}~Remark \ref{remark:cheap}. In such situations,
the operators $\tilde{A}_{\kappa\kappa}$
need not be stored, and the inversion algorithm only requires $O(\sqrt{N}\,\log N)$ memory.
Moreover, if equation (\ref{eq:Ax=b}) is to be solved for several different right hand
sides, subsequent solutions are obtained almost for free via application of the
pre-computed operator $\tilde{A}_{nn}^{-1}$.

Finally we note that while it would require $O(N)$ memory to store $A$ itself,
the scheme only accesses each entry of $A$ once. This means that these elements
can either be computed on the fly (if given by a formula), or read sequentially
from slow memory (``tape'').

\begin{remark}
In cases where each Schur complement $\tilde{A}_{\kappa\kappa}$ is not only
compressible in the ``HSS''-sense, but also in the ``SSS''-sense,
the cost of approximately inverting $\tilde{A}_{\kappa\kappa}$ can be reduced from $O(\kappa^{3})$
to $O(\kappa)$, see Section \ref{sec:prel}. The total cost of Algorithm I is then
$$
T_{\rm total}^{\rm SSS} \sim \sum_{\kappa=1}^{n} \kappa \sim n^{2} \sim N.
$$
This is typically the case when $A$ results from the discretization
of an elliptic partial differential equation.
\end{remark}


\section{Generalizations}

The scheme presented here can in principle be adapted to more
general grids in two and three dimensions. The generalization to
other difference operators on uniform square grids is trivial. Other
two-dimensional grids that are uniform in the sense that they can
easily be partitioned into a sequence of concentric annuli can also
quite easily be handled, and we expect the performance of such
schemes to be similar to the performance reported in Section
\ref{sec:numerics}.

For grids arising from adaptive mesh-refinement, or involve more complex geometries,
it is still possible to construct $O(N \log^{\kappa} N)$ inversion schemes; but they
will very likely require algorithms involving a broader palette of operations on
compressible matrices such as matrix-matrix multiplications. What is interesting
about the specific method given here is that it requires only matrix-inversions
and diagonal updates.

%

\section{Numerical examples}
\label{sec:numerics}

The $O(N\log^{2}N)$ numerical scheme described in Sections \ref{sec:babuska}
and \ref{sec:fast} has been implemented and tested on a conduction problem
on a square uniform grid. We assigned each bar in the grid a conductivity drawn
from a uniform random distribution on the interval $[1,\,2]$. For a range of
grid sizes between $50 \times 50$ and $1\,000 \times 1\,000$, we computed the
operator $\tilde{A}_{n}^{-1}$ described in Section
\ref{sec:babuska}. The computational cost, the amount of memory required, and
the accuracies obtained are presented in Table \ref{table:main}. The following
quantities are reported:
\begin{tabbing}
\mbox{}\hspace{3mm}\=\hspace{12mm}\= \kill
\>$T_{\rm invert}$ \> Time required to construct $\tilde{A}_{n}$ (in seconds)\\
\>$T_{\rm apply}$  \> Time required to apply $\tilde{A}_{\rm n}$ (in seconds)\\
\>$M$              \> Memory required to construct $\tilde{A}_{n}$ (in kilobytes)\\
\>$e_{1}$          \> The largest error in any entry of $\tilde{A}_{n}^{-1}$\\
\>$e_{2}$          \> The error in $l^{2}$-operator norm of $\tilde{A}_{n}^{-1}$\\
\>$e_{3}$          \> The $l^{2}$-error in the vector $\tilde{A}_{nn}^{-1}\,r$ where $r$ is a unit vector
                      of random direction.\\
\>$e_{4}$          \> The $l^{2}$-error in the first column of $\tilde{A}_{nn}^{-1}$.
\end{tabbing}
We estimated $e_{1}$ and $e_{2}$ by comparing the result from the fast algorithm with
the result from a brute force calculation of the Schur complement. We estimated $e_{3}$
and $e_{4}$ by solving equation (\ref{eq:Ax=b}) using iterative methods.

Some technical notes:
\begin{itemize}
\item The experiments were run on a 2.8GHz Pentium 4 PC with 512Mb of RAM.
\item Off-diagonal blocks in the HSS-representations were represented to the
      fixed accuracy $\varepsilon = 10^{-7}$.
\item The ranks used in the off-diagonal blocks was allowed to vary from block
      to block (it was determined adaptively).
\item The code used is written in a Matlab-FORTRAN hybrid. It is not at all
      optimized. Significant gains in efficiency should be obtainable by choosing
      block sizes in more intelligently than we did.
\end{itemize}

\begin{figure}
\begin{tabular}{ccc}
\includegraphics[width=45mm]{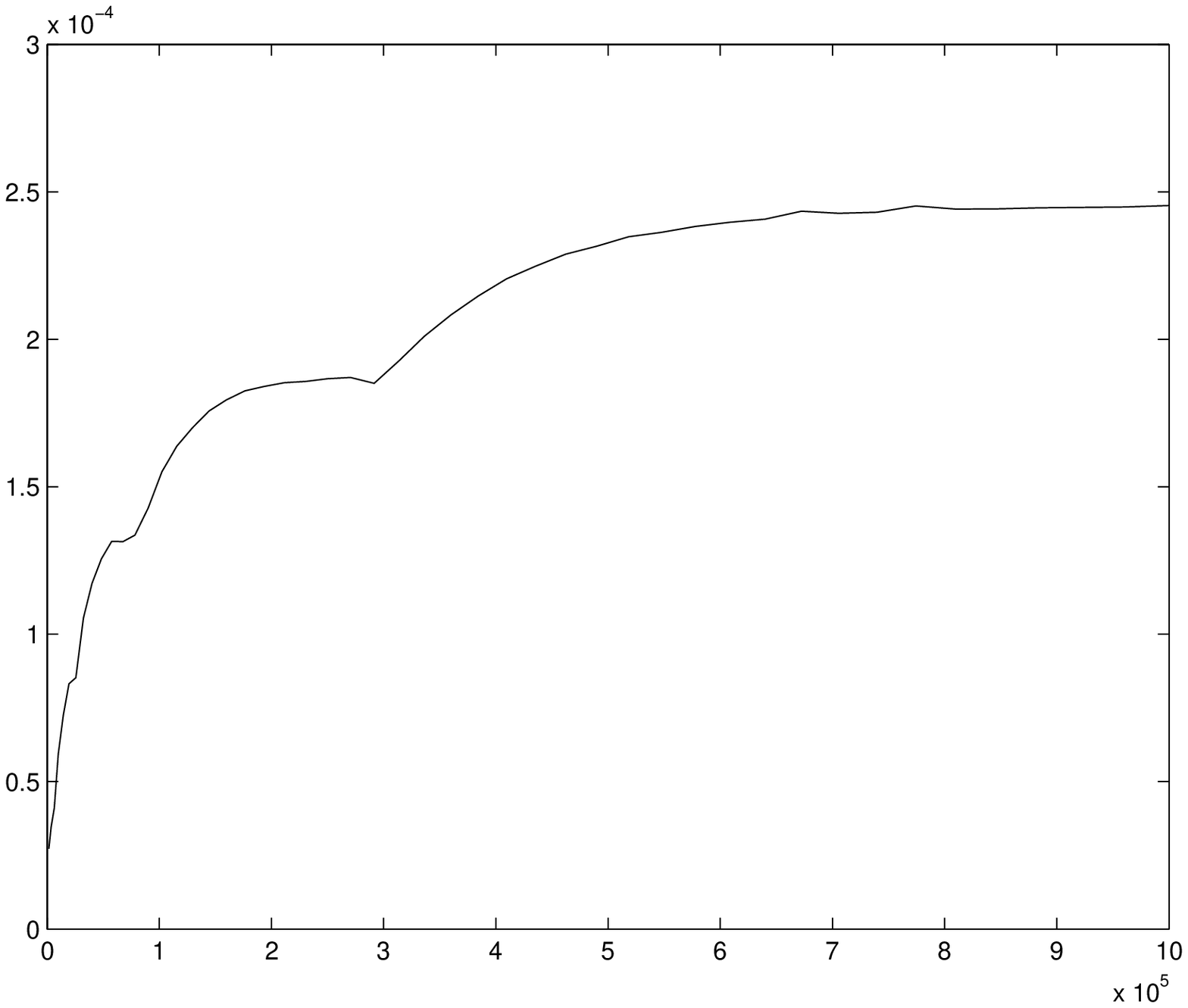}&
\includegraphics[width=45mm]{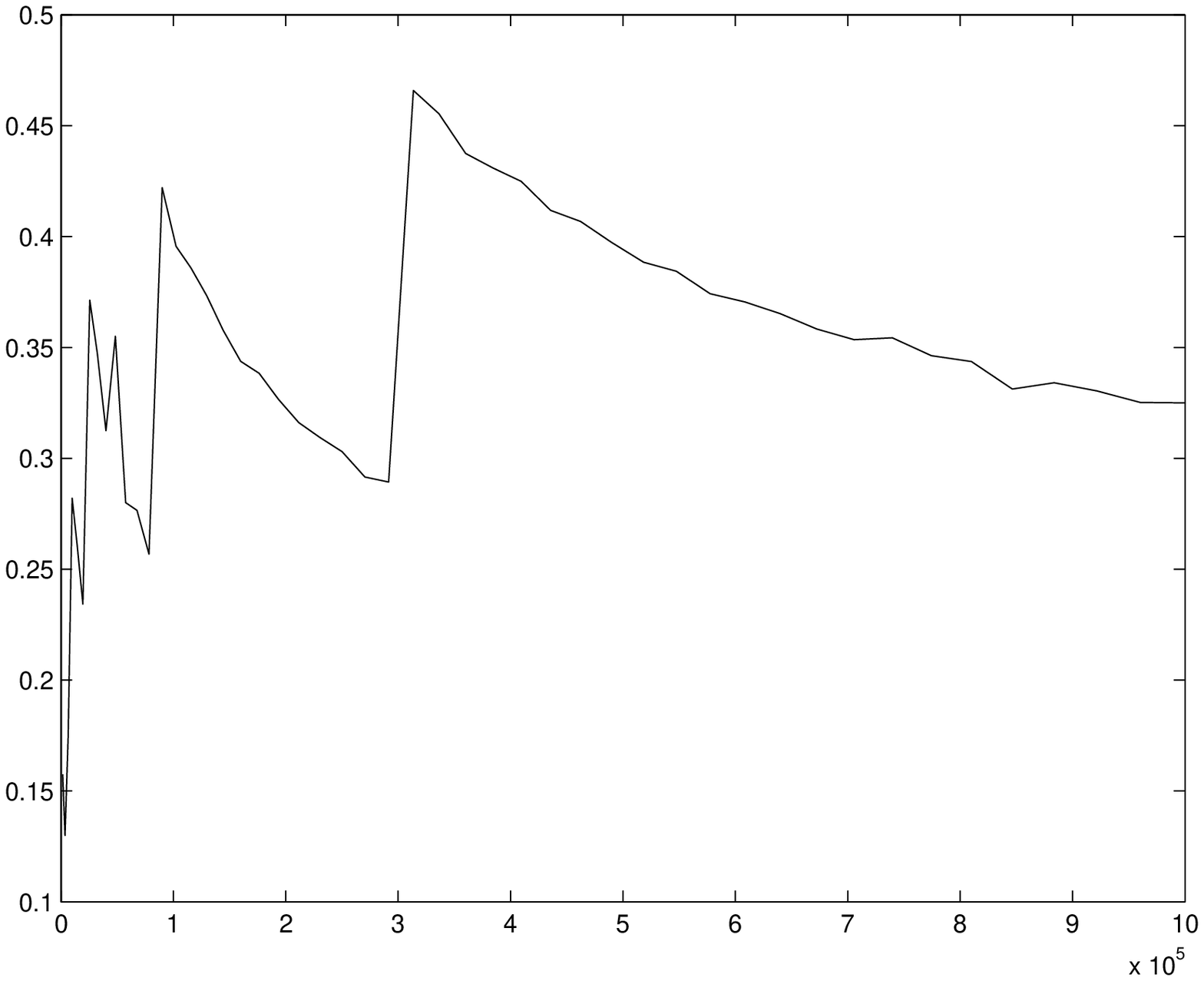}&
\includegraphics[width=45mm]{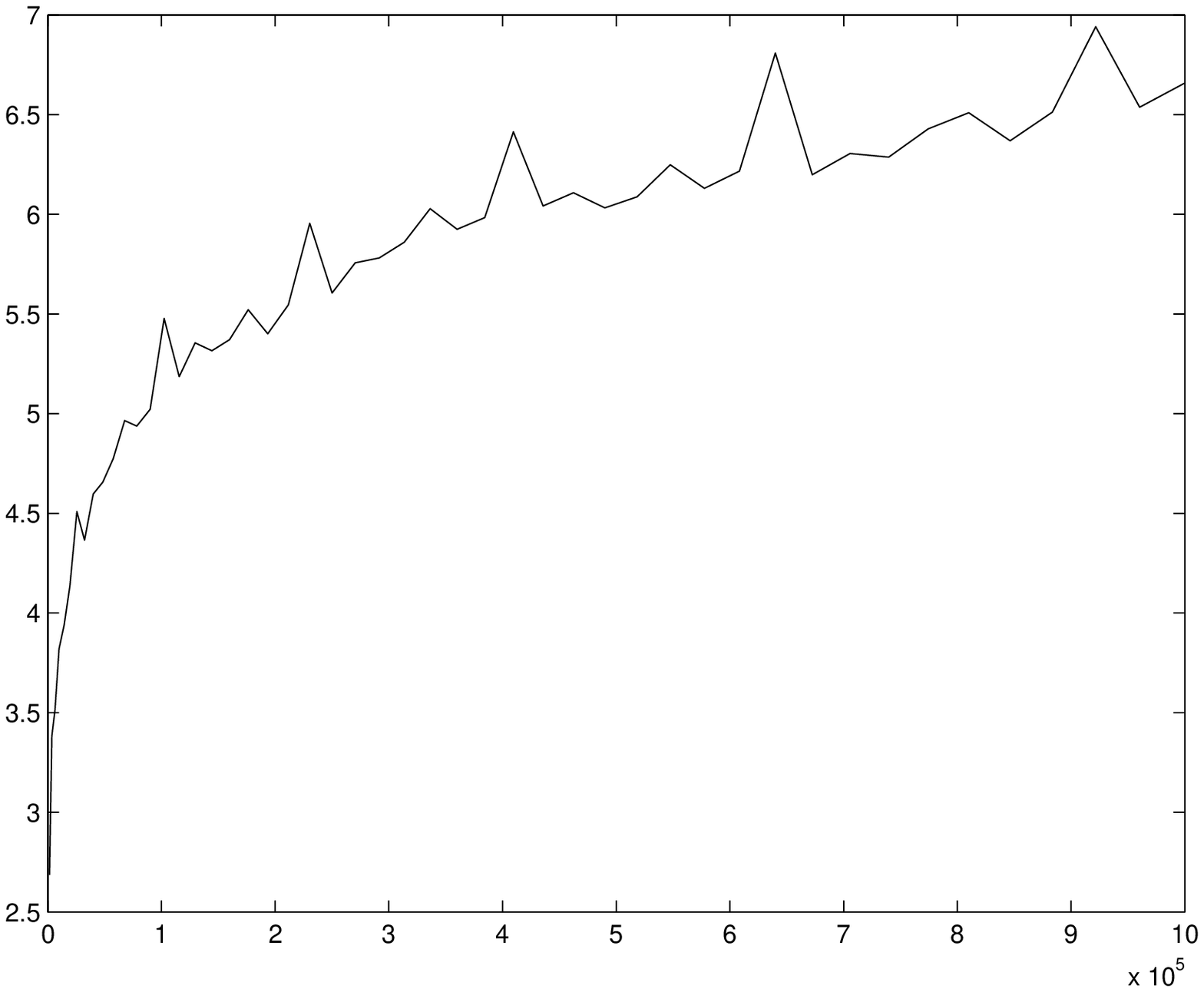}\\
(a) & (b) & (c)
\end{tabular}
\caption{Plots of (a) $T_{\rm invert}/N$ versus $N$, (b) $T_{\rm apply}/\sqrt{N}$ versus $N$,
(c) $M/\sqrt{N}$ versus $N$.}
\label{fig:scaling}
\end{figure}

\begin{figure}
\includegraphics[height=85mm]{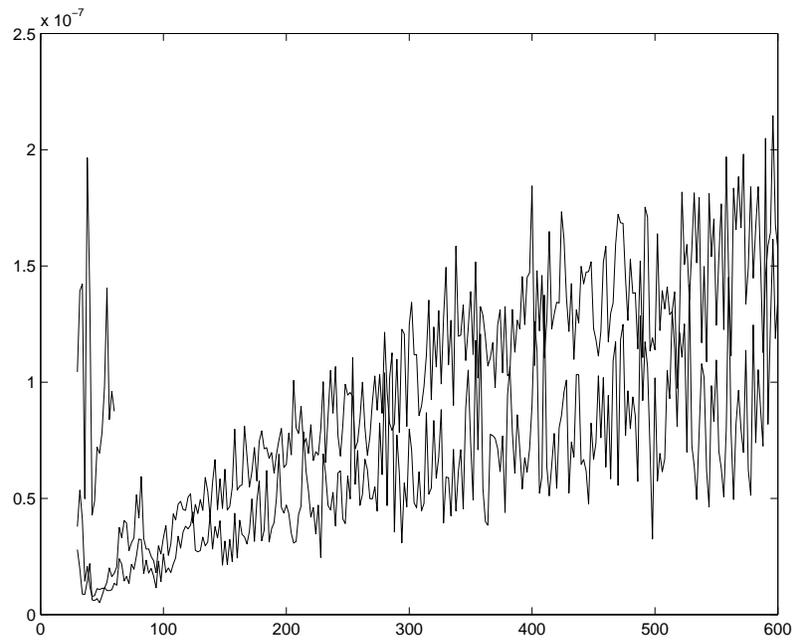}
\caption{Plot of the errors versus \textit{the square root} of the problem size $N$. The short top line on the left
is $e_{2}$, the two lines extending to the right are $e_{3}$ and $e_{4}$, with $e_{4}$
on the top.}
\label{fig:errors}
\end{figure}

The scaling of $T_{\rm invert}$, $T_{\rm apply}$, and $M$ with $N$ is displayed in Figure \ref{fig:scaling}.
These tables appear to support our claims regarding the performance of the scheme.
The values of $e_{1}$, $e_{2}$, $e_{3}$, and $e_{4}$ for different values of $N$ are shown in Figure \ref{fig:errors}.
The table appears to indicate that errors grow as the square root of $N$ (in other words,
linearly with the number of steps performed in Algorithm 1). This means that we could
easily keep the error constant by very moderately increase $\varepsilon$ as the problem
size gets larger.

\begin{figure}
\includegraphics[height=85mm]{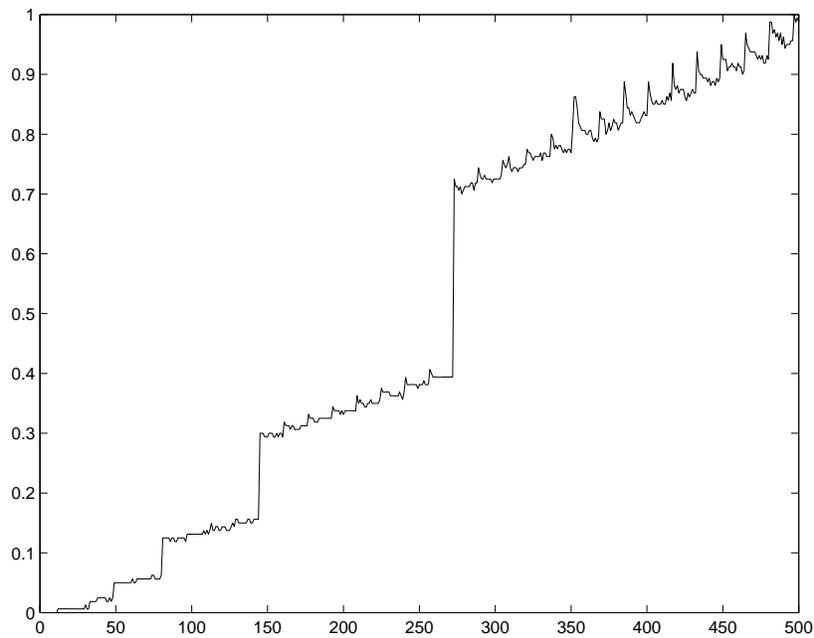}
\caption{Plot of $t_{\kappa}$ versus $\kappa$.}
\label{fig:stepbystep}
\end{figure}

Figure \ref{fig:stepbystep} gives the time $t_{\rm kappa}$ (in seconds) required to
perform step $\kappa$ of the fast version of Algorithm 1. The large jumps in the curve
correspond to repartitioning of the HSS-matrices. (The jaggedness of the curve
gives an indication of how poorly optimized the code is.)

\begin{figure}
\begin{tabular}{r|rrrrrrr}
$N$ & $T_{\rm invert}$ & $T_{\rm apply}$ & $M$ & $e_{1}$ & $e_{2}$ & $e_{3}$ & $e_{4}$ \\ \hline
  10000 & 5.93e-1 & 2.82e-3 & 3.82e+2 & 1.29e-8 & 1.37e-7 & 2.61e-8 & 3.31e-8 \\
  40000 & 4.69e+0 & 6.25e-3 & 9.19e+2 & 9.35e-9 & 8.74e-8 & 4.71e-8 & 6.47e-8 \\
  90000 & 1.28e+1 & 1.27e-2 & 1.51e+3 &     --- &     --- & 7.98e-8 & 1.25e-7 \\
 160000 & 2.87e+1 & 1.38e-2 & 2.15e+3 &     --- &     --- & 9.02e-8 & 1.84e-7 \\
 250000 & 4.67e+1 & 1.52e-2 & 2.80e+3 &     --- &     --- & 1.02e-7 & 1.14e-7 \\
 360000 & 7.50e+1 & 2.62e-2 & 3.55e+3 &     --- &     --- & 1.37e-7 & 1.57e-7 \\
 490000 & 1.13e+2 & 2.78e-2 & 4.22e+3 &     --- &     --- &     --- &     --- \\
 640000 & 1.54e+2 & 2.92e-2 & 5.45e+3 &     --- &     --- &     --- &     --- \\
 810000 & 1.98e+2 & 3.09e-2 & 5.86e+3 &     --- &     --- &     --- &     --- \\
1000000 & 2.45e+2 & 3.25e-2 & 6.66e+3 &     --- &     --- &     --- &     ---
\end{tabular}
\caption{Table summarizing the computational experiment described in Section \ref{sec:numerics}.}
\label{table:main}
\end{figure}

\bibliography{main_bib}
\bibliographystyle{amsplain}

\begin{appendix}

\section{An $O(n\log^{2}n)$ inversion scheme for HSS-matrices}
\label{sec:invert}

A recursive fast inversion scheme for HSS matrices follows
immediately from the following formula for the inverse of a
$2\times 2$ block matrix:
\begin{multline*}
\left[\begin{array}{cc}
A_{11} & A_{12} \\
A_{21} & A_{22}
\end{array}\right]^{-1} =\\
\left[\begin{array}{cc}
\bigl(A_{11} - A_{12}\,A_{22}^{-1}\,A_{21}\bigr)^{-1} &
-\bigl(A_{11} - A_{12}\,A_{22}^{-1}\,A_{21}\bigr)^{-1}\,A_{12}\,A_{22}^{-1} \\
-A_{22}^{-1}\,A_{21}\,\bigl(A_{11} - A_{12}\,A_{22}^{-1}\,A_{21}\bigr)^{-1} &
A_{22}^{-1} + A_{22}^{-1}\,A_{21}\,\bigl(A_{11} - A_{12}\,A_{22}^{-1}\,A_{21}\bigr)^{-1}\,A_{12}\,A_{22}^{-1}
\end{array}\right]
\end{multline*}
From the formula, we immediately get Algorithm 2, displayed in Figure \ref{fig:algorithm2}.
To see that this algorithm has complexity $O(n \log^{2} n)$, we note that the matrices
$A_{12}$ and $A_{21}$ all have low rank. This means that the matrix-matrix multiplications
that occur on lines (7) and (9) in fact consist simply of a small number of multiplications
between HSS-matrices and vectors. Moreover, the matrix additions in lines (7) and (9) are in
fact low-rank updates to HSS-matrices.

\begin{figure}
\begin{center}
\fbox{
\begin{minipage}{120mm}
\begin{center}
\underline{\textit{Algorithm 2:}}
\end{center}

\lsp

\begin{tabbing}
\hspace{10mm} \= \hspace{3mm} \= \hspace{3mm}\= \kill
(1)\> \textbf{function} $B$ = \textit{invert\underline{\mbox{ }}HSS\underline{\mbox{ }}matrix}($A$)\\
\\
(2)\>\> \textbf{if} ($A$ is ``small'') \textbf{then}\\
(3)\>\>\> $B = A^{-1}$\\
(4)\>\>\textbf{else}\\
(5)\>\>\> Split $A = \left[\begin{array}{cc} A_{11} & A_{12} \\ A_{21} & A_{22} \end{array}\right]$.\\
(6)\>\>\> $X_{22}$ = \textit{invert\underline{\mbox{ }}HSS\underline{\mbox{ }}matrix}($A_{22}$)\\
(7)\>\>\> $Y_{11} = A_{11} - A_{12}\,X_{22}\,A_{21}$\\
(8)\>\>\> $X_{11}$ = \textit{invert\underline{\mbox{ }}HSS\underline{\mbox{ }}matrix}($Y_{11}$)\\
(9)\>\>\> $B = \left[\begin{array}{cc} X_{11} & -X_{11}\,A_{12}\,X_{22} \\
                                    -X_{22}\,A_{21}\,X_{11} & X_{22} + X_{22}\,A_{21}\,X_{11}\,A_{12}\,X_{22}\end{array}\right]$.\\
(10)\>\>\textbf{end if}\\
\\
(11)\>\textbf{end function}
\end{tabbing}
\end{minipage}}
\end{center}
\caption{Algorithm for inverting HSS matrices. Note that $A_{11}$, $A_{22}$, $X_{11}$, $X_{22}$, and $Y_{11}$
are all HSS-matrices, and that $A_{12}$ and $A_{21}$ are low-rank matrices.}
\label{fig:algorithm2}
\end{figure}
\end{appendix}


\end{document}